\documentclass[11pt]{article}
\usepackage{amssymb}
\usepackage{amsmath}
\usepackage{amsfonts}
\usepackage{graphicx}
\usepackage{subfig}

\numberwithin{equation}{section}

\begin{document}

\title{{\huge \textbf{Evolution of Spacelike Curves and Special Timelike Ruled
Surfaces in the Minkowski Space }}}
\author{ Dae Won Yoon $^\dag$, Z\"{u}hal K\"{u}\c{c}\"{u}karslan Y\"{u}zba\c{s}\i $^\ddag$ and Ebru Cavlak Aslan$^\S $ }
\date{}
\maketitle
{\footnotesize {%
\centerline  {$\dag$ Department of Mathematics Education and
RINS}}}

{\footnotesize \centerline  {Gyeongsang National University }}

{\footnotesize \centerline{ Jinju 52828, Republic of Korea} }

{\footnotesize \centerline { {E-mail address}:   {\tt dwyoon@gnu.ac.kr}} }
{\footnotesize \vskip 0.2 cm {%
\centerline  { $\ddag$, $\S$ Department of
Mathematics}}} {\footnotesize {\centerline  {F\i rat University }}}

{\footnotesize \centerline{ 23119 Elazig, Turkey} }

{\footnotesize 
\centerline { {E-mail address}:   {\tt
zuhal2387@yahoo.com.tr, ebrucavlak@hotmail.com}} }

\begin{abstract}
In this paper, we get the time evolution equations of the curvature and
torsion of the evolving spacelike curves in the Minkowski space. Also, we
give inextensible evolutions of timelike ruled surfaces that are produced by
the timelike normal and spacelike binormal vector fields of spacelike curve and derive the
necessary conditions for an inelastic surface evolution. Then, we compute
the coefficients of the first and second fundamental forms, the Gauss and
mean curvatures for timelike special ruled surfaces. As a result, we give applications of the evolution equations for the curvatures of the curve in terms of the velocities and get
the exact solutions for these new equations.
\end{abstract}

\renewcommand{\thefootnote}{} \footnote{%
}

\footnote{

The first author was supported by Basic Science Research Program through the National Research Foundation of Korea(NRF) funded by the Ministry of Education (NRF-2018R1D1A1B07046979)}

\renewcommand{\thefootnote}{\arabic{footnote}} \setcounter{footnote}{0}

\section{Introduction}

The time evolution of a curve or a surface is generated by flows, in
particular inextensible flows of a curve or a surface. The flow of a curve
or a surface is said to be inextensible if its arclength is preserved or
the intrinsic curvature is preserved, respectively. Physically, the
inextensible curve flows lead to motions in which no strain energy is
induced. Also, the evolutions of curves have many important applications of
physics as magnetic spin chains and vortex filaments \cite{balak,hasi,lak}.

The problems of how to get the evolution of the curves or surfaces in time is of
deep interest and have been studied in different spaces by many researchers.
Kwon et al. derived the corresponding equations for inextensible flows of
developable surfaces, \cite{kwon}. Hussien et al. obtained the evolutions of
the special surfaces rely on the evolutions of their directrices, \cite{huss}. Recently,
many authors have studied geometric flow problems on the curves or surfaces 
\cite{abdel, abdel2,alk,bektas,gurbuz,lamb,yil,kucuk}.

In this study, we get the evolution of curves via the velocities of the
moving frame in Minkowski space. We also classify the special timelike ruled
surfaces on the evolving spacelike curve where the generator is the timelike
principal and spacelike binormal vectors to the spacelike curves on the
timelike surfaces. We give the necessary conditions for the inelastic special
ruled surface evolutions and compute the gauss and mean curvatures for
them. Then, we obtain a pair of coupled nonlinear partial differential
equations governing the time evolution of the curvatures of the evolving
spacelike curve in Minkowski space. We give the new geometric models of the
evolution equations for curvatures from the main equation in Minkowski space and
get the exact solutions for them and show the moving curve for these
solutions. From these exact solutions, we derive two types of nonlinear
traveling solitary wave, which are kink and bell-shaped solitary waves. The
kink solitary wave appears balance between nonlinearity and dissipation
supports known the nonlinear wave of stable shape. The bell-shaped solitary
wave appears in consequence of the balance between nonlinearity and
dispersion. Also, the bell-shaped wave improve to evaluate the dynamic
modulus \cite{choi, biswas}.

\bigskip

\section{Preliminaries}

The 3-dimensional Minkowski space $E_{1}^{3}$ is the real vector space $%
E^{3}$ provided with the Lorentzian inner product given%
\begin{equation*}
\left\langle a,b\right\rangle =a_{1}b_{1}+a_{2}b_{2}-a_{3}b_{3},
\end{equation*}%
with $a=(a_{1},a_{2},a_{3})$ and $b=(b_{1},b_{2},b_{3})$.

A vector $a\in E_{1}^{3}$ is said to be a spacelike vector when $%
\left\langle a,a\right\rangle >0$ or $a=0$. It is said to be a timelike and
a null (light-like) vector in case of $\left\langle a,a\right\rangle <0$, and $%
\left\langle a,a\right\rangle =0$ for $a\neq 0$, respectively. Similarly, a
curve is said to be spacelike, if its velocity vector is spacelike. For
furthermore information, we refer to \cite{oneil}.

Let $r =r (s)$ be a spacelike curve in 3-dimensional Minkowski
space $E_{1}^{3}$. Then Frenet formulas are given by%
\begin{equation}
\left\{ 
\begin{array}{c}
T_{s}=\kappa N \\ 
N_{s}=\kappa T+\tau B \\ 
B_{s}=\tau N%
\end{array}%
\right. ,  \label{1}
\end{equation}%
where $T=r _{s},$ $N$ and $B$ are called the vectors of the unit
spacelike tangent, timelike principal normal and spacelike binormal of $%
r(s)$, respectively. Also $\kappa $ and $\tau $ are geometric
parameters that represent, respectively, the curvature and torsion of $%
r(s),$ \cite{woe}. Through this paper, the subscripts describe partial
derivatives.

We also know that a curve is uniquely given by two scalar quantities, so-called
curvature and torsion.

If $r(s)$ moves with time $t,$ then \eqref{1} is becoming functions of
both $s$ and $t$. We can give the evolution equations of $\{T,N,B\}$ quite
generally, in a form similar to \eqref{1} as following \cite{ugurlu}

\begin{equation}
\left\{ 
\begin{array}{c}
T_{t}=\alpha N-\beta B \\ 
N_{t}=\alpha T+\gamma B \\ 
B_{t}=\beta T+\gamma N%
\end{array}%
\right. .  \label{2}
\end{equation}%
Clearly, $\alpha ,\beta $ and $\gamma $ (which are the velocities of the
moving frame) detect the motion of the $\alpha $.

Let $x=x(s,t)$ be a surface parametrization in $E_{1}^{3}$. Then, the vectors $%
x_{s}$ and $x_{t}$ are tangential to $M$ at $p$. Let $U$ be the standard unit 
normal vector field on a surface defined by 
\begin{equation}
U=\frac{x_{s}\wedge x_{t}}{\left\Vert x_{s}\wedge x_{t}\right\Vert }.
\label{2b}
\end{equation}%
Then the first and the second fundamental forms of a surface are given by 
\cite{carmo, lopez}, respectively,%
\begin{eqnarray*}
I &=&Eds^{2}+2Fdsdt+Gdt^{2}, \\
II &=&eds^{2}+2fdsdt+gdt^{2},
\end{eqnarray*}%
where%
\begin{eqnarray}
E &=&\left\langle x_{s},x_{s}\right\rangle ,\text{ }F=\left\langle
x_{s},x_{t}\right\rangle ,\text{ }G=\left\langle x_{t},x_{t}\right\rangle ,
\label{2a} \\
e &=&\left\langle x_{ss},U\right\rangle ,\text{ }f=\left\langle
x_{st},U\right\rangle ,\text{ }g=\left\langle x_{tt},U\right\rangle .  \notag
\end{eqnarray}%
Denote $W=EG-F^{2}$. The surface is spacelike if $W>0$ and it is timelike if 
$W<0$. We give the notation $\varepsilon =\left\langle U,U\right\rangle =\pm
1 $. Therefore, we can write%
\begin{equation*}
\left\Vert x_{s}\wedge x_{t}\right\Vert =\sqrt{-\varepsilon W}.
\end{equation*}

Thus, the Gauss and the mean curvatures are defined by%
\begin{eqnarray}
K &=&\varepsilon \frac{eg-f^{2}}{EG-F^{2}},  \label{2c} \\
H &=&\varepsilon \frac{Eg-2Ff+Ge}{2\left( EG-F^{2}\right) }.  \notag
\end{eqnarray}

On the other hand, a surface evolution $x(u,s,t)$ and its flow $\frac{%
\partial x}{\partial t}$ are said to be inextensible if%
\begin{equation}
\frac{\partial E}{\partial t}=\frac{\partial F}{\partial t}=\frac{\partial G%
}{\partial t}=0.  \label{2d}
\end{equation}

\section{Time Evolution Equations of the Spacelike Curves and Special
Timelike Ruled Surfaces in $E_{1}^{3}$}

For spacelike inextensible curves, the moving frame must be holded the
compatibility conditions%
\begin{equation}
T_{st}=T_{ts},\text{ }N_{st}=N_{ts}\text{ and }B_{st}=B_{ts}.  \label{3}
\end{equation}%
Here spacelike inextensible curves mean that the flow described by the
equations \eqref{2} preserves the curves in arc-length parametrization.

If we substitute \eqref{1} and \eqref{2} into \eqref{3}, then we get%
\begin{eqnarray*}
\kappa \alpha T+\kappa _{t}N+\gamma \kappa B &=&\alpha \kappa T+\left(
\alpha _{s}-\tau \beta \right) N+\left( -\beta _{s}+\tau \alpha \right) B, \\
\left( \kappa _{t}+\tau \beta \right) T+\left( \kappa \alpha +\tau \gamma
\right) N+\left( \tau _{t}-\kappa \beta \right) B &=&\alpha _{s}T+\left(
\kappa \alpha +\tau \gamma \right) N+\gamma _{s}B, \\
\kappa \tau T+\tau _{t}N+\tau \gamma B &=&\left( \beta _{s}+\kappa \gamma
\right) T+\left( \kappa \beta +\gamma _{s}\right) N+\tau \gamma B.
\end{eqnarray*}%
From the above equations we get%
\begin{eqnarray}
\kappa _{t} &=&\alpha _{s}-\tau \beta ,  \label{4} \\
\tau _{t} &=&\gamma _{s}+\kappa \beta ,  \notag \\
\gamma &=&\frac{\alpha \tau -\beta _{s}}{\kappa }.  \notag
\end{eqnarray}%
The temporal evolution of $\kappa $ and $\tau $ of a spacelike curve can be
given in terms of $\{\alpha ,\beta ,\gamma \}$ which is obtained as
\begin{eqnarray}
\kappa _{t} &=&\alpha _{s}-\tau \beta ,  \label{5} \\
\tau _{t} &=&\left( \frac{\alpha \tau -\beta _{s}}{\kappa }\right)
_{s}+\kappa \beta .  \notag
\end{eqnarray}%
The equation \eqref{5} is one of the main results of this paper. We determine the
equations of motion of the spacelike curve for a given $\{\alpha ,\beta
,\gamma \}.$ Then, we choose $\{\alpha ,\beta ,\gamma \}$ in terms of the $%
\left\{ \kappa ,\tau \right\} .$

\section{Inextensible Flows of the Timelike Normal and Binormal Ruled
Surfaces}

\subsection{Timelike Normal Surfaces}

We can give a timelike normal surface $M$ in $E_{1}^{3}$ with the parametric
representation in a similar way to \cite{huss} as follows 
\begin{equation*}
x(s,u,t)=r (s,t)+uN(s,t),
\end{equation*}%
where a spacelike curve $r$ $=r(s)$ and a timelike principal
normal $N=N(s)$ of $r$ move with time $t$ are showed $r(s,t)$
and $N(s,t).$

We get the first and second partial derivatives of $M$ in terms of $s$ and $%
u $ are computed as follows 
\begin{eqnarray}
x_{s}(s,u,t) &=&(1+u\kappa )T+\tau uB,  \label{6} \\
x_{u}(s,u,t) &=&N(s,t),  \notag \\
x_{ss}(s,u,t) &=&u\kappa _{s}T+\left( (1+u\kappa )\kappa +\tau ^{2}u\right)
N+\tau _{s}uB,  \notag \\
x_{su}(s,u,t) &=&\kappa T+\tau B,  \notag \\
x_{uu}(s,u,t) &=&0.  \notag
\end{eqnarray}%
Also, using \eqref{2b} and cross product in $E_{1}^{3},$ the unit normal
vector field on $M$ is found as%
\begin{equation*}
U=\frac{\tau uT+(1+u\kappa )B}{\sqrt{\tau ^{2}u^{2}+(1+u\kappa )^{2}}}.
\end{equation*}%
Then, from \eqref{2a}, the coefficients of the first and second fundamental
forms can be expressed as 
\begin{eqnarray*}
E &=&\tau ^{2}u^{2}+(1+u\kappa )^{2},F=0,G=-1, \\
e &=&\frac{\tau u^{2}\kappa _{s}+\tau _{s}u(1+u\kappa )}{\sqrt{\tau
^{2}u^{2}+(1+u\kappa )^{2}}},f=\frac{2\tau u\kappa +\tau }{\sqrt{\tau
^{2}u^{2}+(1+u\kappa )^{2}}},g=0.
\end{eqnarray*}%
Thus, we can compute $K$ and $H$ using \eqref{2c} as follows%
\begin{eqnarray}
K &=&\frac{\left( 2\tau u\kappa +\tau \right) ^{2}}{\left( \tau
^{2}u^{2}+(1+u\kappa )^{2}\right) ^{2}},  \label{6b} \\
H &=&\frac{\tau u^{2}\kappa _{s}+\tau _{s}u(1+u\kappa )}{2\left( \tau
^{2}u^{2}+(1+u\kappa )^{2}\right) ^{\frac{3}{2}}}.  \notag
\end{eqnarray}%
With the help of \eqref{2d}, if the normal surface is inextensible, then one
can derive 
\begin{equation}
\tau \tau _{t}u+(1+u\kappa )\kappa _{t}=0.  \label{6c}
\end{equation}

\subsection{Timelike Binormal Surfaces}

We can give a timelike binormal surface in $E_{1}^{3}$ with the parametric
representation in a similar way to \cite{huss} as follows 
\begin{equation*}
x(s,v,t)=r(s,t)+vB(s,t),
\end{equation*}%
where a spacelike curve $r $ $=r(s)$ and a spacelike principal
binormal $B=B(s)$ of $r$ move with time $t$ are showed $r(s,t)$
and $B(s,t).$ By following a similar way as above, the coefficients of the
first and second fundamental forms are obtained as 
\begin{eqnarray}
E &=&1-v^{2}\tau ^{2}, F=0, G=1,  \label{6d} \\
e &=&-\frac{\tau ^{2}v^{2}\kappa +(\kappa +v\tau _{s})}{\sqrt{v^{2}\tau
^{2}-1}},f=-\frac{\tau }{\sqrt{v^{2}\tau ^{2}-1}},g=0.  \notag
\end{eqnarray}%
Thus, we can compute $K$ and $H$  using \eqref{2c} as follows%
\begin{eqnarray}
K &=&\frac{\tau ^{2}}{\left( v^{2}\tau ^{2}-1\right) ^{2}},  \label{6e} \\
H &=&\frac{\tau ^{2}v^{2}\kappa +(\kappa +v\tau _{s})}{2\left( v^{2}\tau
^{2}-1\right) ^{\frac{3}{2}}}.  \notag
\end{eqnarray}%
By taking account of \eqref{2d}, if the binormal surface is inextensible,
then one can derive 
\begin{equation}
v^{2}\tau \tau _{t}=0.  \label{6f}
\end{equation}

\section{Applications}

\textbf{Type 1} The evolution equations for the curvature and the torsion of
the curve in terms of the velocities $\left\{ \alpha ,\beta ,\gamma \right\}
=\left\{ 0,\kappa ,-\dfrac{\kappa _{s}}{\kappa }\right\} $ obtained using %
\eqref{5} as%
\begin{eqnarray}
\kappa _{t} &=&-\tau \kappa ,  \label{7} \\
\tau _{t} &=&\left( \frac{\kappa _{s}^{2}-\kappa \kappa _{ss}}{\kappa ^{2}}%
\right) +\kappa ^{2}.  \notag
\end{eqnarray}%
We seek solutions of \eqref{7} in the form%
\begin{eqnarray}
\kappa (s,t) &=&A_{1}\tanh ^{p_{1}}\xi ,  \label{8} \\
\tau (s,t) &=&A_{2}\tanh ^{p_{2}}\xi ,  \notag
\end{eqnarray}%
where $\xi =\eta (s-\upsilon t)$ . In \eqref{8}, $%
A_{1},A_{2}$ and $\eta $ are arbitrary real constants and $\upsilon $ is the
velocity of the solitary wave, we refer to \cite{Ozkan, Inc} for a more information of these solution methods. So, from \eqref{7} and \eqref{8}, we have
\begin{figure}[tbh]
{\ \centering
\includegraphics[scale=0.5]{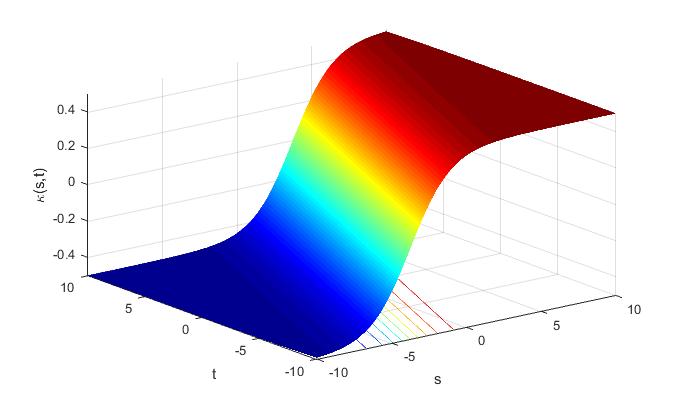} \hspace*{1cm} \centering
}
\caption{The kink soliton solution of $\protect\kappa(s,t)$ for $A_{1}=0.5,
A_{2}=1.$}
\end{figure}
\begin{figure}[tbh]
{\ \centering
\includegraphics[scale=0.5]{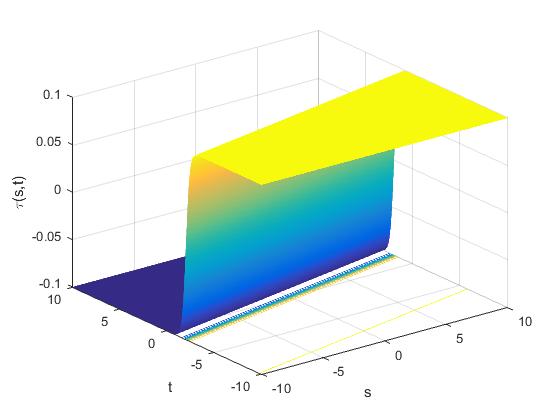} \hspace*{1cm} \centering
}
\caption{The kink soliton solution of $\protect\tau(s,t)$ for $A_{1}=0.1,
A_{2}=0.1.$}
\end{figure}

\begin{equation}
-A_{1}\eta \upsilon p_{1}\tanh ^{p_{1}-1}\xi -A_{1}\eta \upsilon p_{1}\tanh
^{p_{1}+1}\xi +A_{1}A_{2}\tanh ^{p_{1}+p_{2}}\xi =0  \label{9}
\end{equation}%
and 
\begin{equation*}
A_{1}^{2}\eta ^{2}p_{1}^{2}\tanh ^{p_{1}+1}\xi (1+\tanh ^{2}\xi
)^{2}-A_{1}^{2}\eta ^{2}p_{1}(p_{1}-1)\tanh ^{2p_{1}-2}\xi (1+\tanh ^{2}\xi
)+A_{1}^{4}\tanh ^{4p_{1}}\xi
\end{equation*}%
\begin{equation}
-A_{1}^{2}\eta ^{2}p_{1}(p_{1}+1)\tanh ^{2p_{1}}\xi (1+\tanh ^{2}\xi
)+A_{1}^{2}A_{2}\eta \upsilon p_{2}\tanh ^{2p_{1}+p_{2}-1}\xi (1+\tanh
^{2}\xi )=0.  \label{10}
\end{equation}%
From \eqref{9} and \eqref{10}, equating the coefficients of $%
(p_{1}+1,2p_{1}) $ and $(2p_{1},2p_{1}+p_{2}-1)$ gives%
\begin{equation}
p_{1}=p_{2}=1.  \label{11}
\end{equation}%
As a result, we obtain%
\begin{equation*}
\upsilon =\frac{\eta }{A_{2}},
\end{equation*}%
where $\eta ^{2}+A_{1}^{2}=1$. 
Hence, we get the kink soliton solution for \eqref{7} as 
\begin{eqnarray}
\kappa (s,t) &=&A_{1}\tanh (\eta (s-\upsilon t)),  \label{12} \\
\tau (s,t) &=&A_{2}\tanh (\eta (s-\upsilon t)).  \notag
\end{eqnarray}

\textbf{Type 2} The evolution equations for the curvatures of the curve in
terms of the velocities $\left\{ \alpha ,\beta ,\gamma \right\} =\left\{
\tau ,\kappa _{s},\dfrac{\tau -\kappa _{ss}}{\kappa }\right\} $ obtained
using \eqref{5} as%
\begin{eqnarray}
\kappa _{t} &=&\tau _{s}-\tau \kappa _{s},  \label{13} \\
\tau _{t} &=&\left( \frac{(2\tau \tau _{s}-\kappa _{sss})\kappa -(\tau
^{2}-\kappa _{ss})\kappa _{s}}{\kappa ^{2}}\right) +\kappa \kappa _{s}. 
\notag
\end{eqnarray}%
\begin{figure}[tbh]
{\ \centering
\includegraphics[scale=0.5]{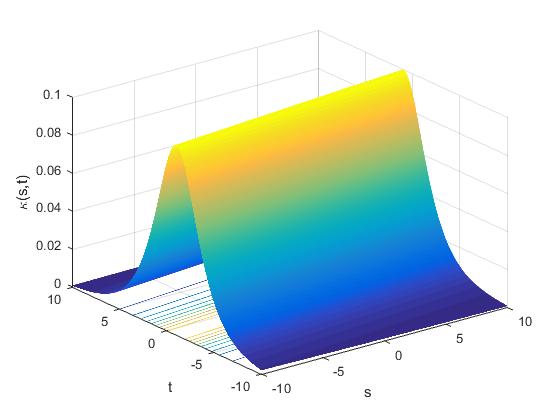} \hspace*{1cm} \centering
}
\caption{The bell-shaped soliton solution of $\protect\kappa(s,t)$ for $%
B_{1}=0.1, B_{2}=-1.$}
\end{figure}
\begin{figure}[tbh]
{\ \centering
\includegraphics[scale=0.5]{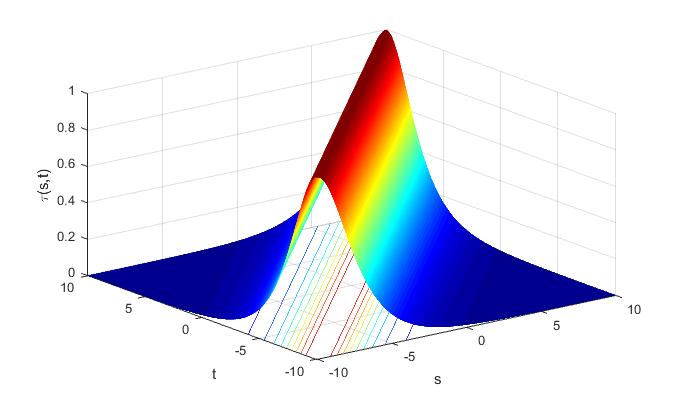} \hspace*{1cm} \centering
}
\caption{The bell-shaped soliton solution of $\protect\tau(s,t)$ for $%
B_{1}=-1, B_{2}=1.$}
\end{figure}
We assume that the solutions are in the form%
\begin{eqnarray}
\kappa (s,t) &=&\frac{B_{1}}{\cosh ^{1/p_{1}}\xi },  \label{14} \\
\tau (s,t) &=&\frac{B_{2}}{\cosh ^{1/p_{2}}\xi },  \notag
\end{eqnarray}%
with $\xi =\eta (s-\upsilon t)$, \cite{biswas}. Substituting \eqref{14} into %
\eqref{13} yields%
\begin{equation}
B_{1}\eta \upsilon \frac{1}{p_{1}}\frac{1}{\cosh ^{1/p_{1}}\xi }+B_{2}\eta 
\frac{1}{p_{2}}\frac{1}{\cosh ^{1/p_{2}}\xi }-B_{1}B_{2}\frac{1}{p_{1}}\frac{%
1}{\cosh ^{1/p_{1}+1/p_{2}}\xi }=0.  \label{15}
\end{equation}%
and 
\begin{equation*}
B_{1}^{2}\eta \frac{1}{p_{1}}\frac{1}{\cosh ^{2/p_{1}}\xi }+B_{2}\eta
\upsilon \frac{1}{p_{2}}\frac{1}{\cosh ^{1/p_{2}}\xi }+\frac{B_{2}^{2}}{B_{1}%
}\eta \left( -\frac{1}{p_{1}}+\frac{2}{p_{2}}\right) \frac{1}{\cosh
^{1/p_{1}+2/p_{2}}\xi }.
\end{equation*}%
\begin{equation}
+2\eta ^{3}\left( -1-\frac{1}{p_{1}}\right) \frac{1}{p_{1}}\frac{1}{\cosh
^{2}\xi }=0.  \label{16}
\end{equation}%
Now, from \eqref{16}, equating the exponent $\left( \dfrac{2}{p_{1}}%
,2\right) $ and $\left( \dfrac{-1}{p_{1}}+\dfrac{2}{p_{2}},\dfrac{1}{p_{2}}%
\right) $ leads to 
\begin{equation}
p_{1}=p_{2}=1.  \label{17}
\end{equation}%
Consequently, we obtain%
\begin{equation}
\upsilon =-\frac{B_{2}}{B_{1}},  \label{18}
\end{equation}%
where $\eta =\dfrac{B_{1}}{2}.$

We obtain following bell-shaped soliton solutions%
\begin{eqnarray}
\kappa (s,t) &=&\frac{B_{1}}{\cosh \xi }=B_{1}\sec h\xi ,  \label{19} \\
\tau (s,t) &=&\frac{B_{2}}{\cosh \xi }=B_{2}\sec h\xi .  \notag
\end{eqnarray}

\section{Conclusion}

In this work, the evolutions of inextensible spacelike curves and timelike special ruled surfaces have been obtained. Then, we derived  the evolution equations for the curvatures of the curve in terms of the velocities and found exact solutions.
We have also found the kink solitary and bell-shaped wave solutions by
ansatz method for the solutions of these equations. In Figs.1-2, the kink solitary waves of $\kappa (s,t)~$ and $%
\tau (s,t)~$ obtained by Eq. \eqref{12} are presented for values $%
A_{1}=0.5, A_{2}=1$ and $A_{1}=0.1, A_{2}=0.1$. The bell-shaped solitary waves
of $\kappa (s,t)~$ and $\tau (s,t)~$ obtained by Eq. \eqref{19} are given by
Figs. 3-4 for $B_{1}=0.1, B_{2}=-1$ and $B_{1}=-1, B_{2}=1$.

\end{document}